\documentclass[12pt]{amsart}
\usepackage{amssymb}

\usepackage{euscript}

\let\cal\mathcal
\input xypic
\xyoption {all}

\makeatletter \@mparswitchfalse \makeatother

\newtheorem{theorem}{Theorem}
\newtheorem{lemma}{Lemma}
\newtheorem{corollary}{Corollary}
\newtheorem{proposition}{Proposition}
\newtheorem{remark}{Remark}

\newtheorem{definition}{Definition}

\def\eqref#1{(\ref{eq#1})}

\def\N{\mathbb N}

\def\M{\cal{M}}
\let\<\langle

\numberwithin{equation}{section}

\begin{document}


\let\\\cr
\let\union\bigcup
\let\inter\bigcap
\def\supp{\operatorname{supp}}
\def\sup{\operatorname{sup}}
\def\inf{\operatorname{inf}}
\def\Im{\operatorname{Im}}
\def\dim{\operatorname{dim}}
\def\Span{\operatorname{span}}
\def\cord{\operatorname{cord}}
\def\lim{\operatorname{lim}}
\def\Re{\operatorname{Re}}
\def\sqn{\operatorname{sqn}}
\def\log{\operatorname{log}}
\def\max{\operatorname{max}}
\def\min{\operatorname{min}}
\let\emptyset\varnothing
\def\exp{\operatorname{exp}}

\title{FACTORIZATION  OF OPERATORS ON
$C^*$-ALGEBRAS}
\author{Narcisse Randrianantoanina}
\address{Department of Mathematics and Statistics, Miami University, Oxford,
Ohio 45056}
\email{randrin@muohio.edu}


\begin{abstract}
Let $\cal A$ be a $C^*$-algebra.  We prove that every
 absolutely summing operator from $\cal A$ into $\ell_2$
factors through a Hilbert space operator that belongs to the
 4-Schatten-von Neumann class.
 We also provide finite dimensional examples that show
that one can not replace the 4-Schatten-von Neumann class by
$p$-Schatten-von Neumann class for any $p <4$.
 As an application, we show that
there exists a modulus of capacity $\varepsilon \to N(\varepsilon)$
 so that if $\cal A$ is a $C^*$-algebra and
 $T \in \Pi_1 (\cal A, \ell_2)$ with $\pi_1 (T) \leq 1$, then
 for every $\varepsilon >0$, the
$\varepsilon$-capacity of the image of the unit ball of $\cal A$ under $T$
does not exceed $N(\varepsilon)$.  This answers positively
 a question raised by Pe\l czy\'nski.
\end{abstract}

\maketitle

\section{Introduction}
It is a well known consequence of a classical result of Grothendieck that
if $X$ is a Banach space and $X^{**}$ is isomorphic to a quotient of a
$C(K)$-space then every absolutely summing operator from $X$ into $\ell_2$
factors through a Hilbert-Schmidt operator.  The present paper is an
attempt to get a generalization of this fact
for the setting of arbitrary $C^*$-algebras.
Different structures  of operators defined on arbitrary $C^*$-algebras was
considered by Pisier in \cite{Pis3} and \cite{Pis4};
 for instance he proved that every
$(p,q)$-summing operators on an arbitrary $C^*$-algebra admit a
factorization similar to that of operators on $C(K)$-spaces,
every operator from any $C^*$-algebra  into any Banach space of cotype 2
factors through Hilbert space.
Using the notion of $C^*$-summing operators introduced by Pisier in
\cite{Pis3}, the author proved in \cite{Ran7},
that absolutely summing operators from $C^*$-algebras into
reflexive spaces are compact.  The main result of this paper states that
for the case of $C^*$-algebras and the range
space  being a
Hilbert space, one can factor every absolutely summing
operator through a Hilbert space operator that belongs to the
 4-Schatten-von Neumann class (see definition below).
The basic idea of the proof of this result is
the factorization of $C^*$-summing operators  used in
\cite{Ran7} and some well known co\"incidence of different classes of
Hilbert space operators.
 This result allows to prove a
quantitative result on the compactness of absolutely summing operators from
$C^*$-algebras into Hilbert spaces, answering a question raised by
Pe\l czy\'nski in \cite{HN} (Problem 3') for the space of compact operators on
Hilbert space. A finite dimensional approach shows that unlike the
commutative case of $C(K)$-spaces, one cannot expect to factor every
absolutely summing operators from  general non-commutative $C^*$-algebras
into Hilbert spaces  through Hilbert-Schmidt operators. In fact our
examples show that the result stated above cannot be improved to the case
of $p$-Schatten-von Neumann class for any $p< 4$.

Our terminology and notation are standard. We refer to \cite{D1} and
\cite{WO} for definitions from Banach space theory and \cite{KR} and
\cite{TAK} for basic properties from $C^*$-algebra and operator algebra
theory.

\section{Preliminaries}

In this section we recall some definitions and facts which we use in the
sequel.
Throughout, the word operator will always mean linear bounded operator
and
$\cal{L}(E, F)$ will stand for the space of all operators from E into F.

\begin{definition}
Let E and F be Banach spaces and $1\leq p < \infty$.  An operator
$T \in \cal L(E,F)$ is said to be absolutely $p$-summing
 (or simply $p$-summing)
if there exists a constant $C$ such that for any finite sequence $(e_1,
e_2, \dots, e_n)$ of $E$, one has
$$\left(\sum^n_{i=1}\Vert T e_i\Vert^p\right)^{\frac 1p} \leq C \sup
\left\{\left(\sum_{i=1}^n \vert \langle e_i,e^*\rangle\vert^p
\right)^{\frac 1p}\ ;\ e^* \in E^*,\ \Vert e^* \Vert \leq
1\right\}.
$$
\end{definition}

The least constant $C$ for the inequality above to hold will be denoted by
$\pi_p (T)$. It is well known that the class of all absolutely $p$-summing
operators from $E$ to $F$ is a Banach space under the norm $\pi_p (.)$.
  This Banach space will be denoted by $\Pi_p (E, F)$ .

\begin{definition}
Let $1\leq q \leq p < \infty$.  An operator $T \in \cal{L}(E, F)$ is said
to be $(p,q)$-summing if there is a constant  $K \geq 0$
for which
$$
\left(\sum^n_{k=1}\Vert T e_i\Vert^p\right)^{\frac 1p} \leq K \sup
\left\{ \left(\sum^n_{i=1} \vert \langle
e^*,e_i\rangle\vert^q\right)^{\frac 1q};\
\ e^* \in E^*,\ \Vert e^*\Vert\leq1\right\}
$$
for every finite sequnce $(e_1, e_2, \dots, e_n)$ in $E$.
\end{definition}

As above, the least constant $K$ for which the inequality holds is the
$(p,q)$-summing norm of $T$ and is denoted by $\pi_{p,q} (T)$.
 The class of
$(p,q)$-summing operators from $E$ into $F$ is a Banach space under the
norm $\pi_{p,q}(.)$. This class will be denoted by $\Pi_{p,q}(E,F)$.

Another class of operators relevant for our discussion is the
 Schatten-von Neumann class.

\begin{definition}
For $1\leq p < \infty$, $H_1$ and $H_2$ Hilbert spaces, the $p$-th
 Schatten-von Neumann class consists of all compact operators $U: H_1
\longrightarrow H_2$ that has a representation of the form
$$
(*)\ \ \ \ \ U = \sum_{n=1}^\infty \alpha_n (\ .\ , e_n) f_n,
$$
where $(e_n)_n$ is an orthonormal sequence in $H_1$, $(f_n)_n$ is an
orthonormal sequence in $H_2$, and $(\alpha_n)_n \in \ell_p$.
\end{definition}

We will refer to $(*)$ as an orthonormal representation of $U$.
 It is well known
that one can always choose the sequence $(\alpha_n)_n$ in the
representation $(*)$ to satisfy $0 \leq \alpha_{n+1} \leq \alpha_n$ for all
admissible indices.  The $p$-th Schatten-von Neumann norm is defined by
$$
\sigma_p(U) = \left(\sum_{n=1}^\infty \vert \alpha_n \vert^p\right)^\frac 1p$$ and the p-th
Schatten-von Neumann class is denoted by $S_p(H_1, H_2)$.

\begin{definition}
Let $E$ and $F$ be Banach spaces, $1\leq p \leq \infty$. We say that an
operator $T \in  \cal L(E,F)$ is $L_p$-factorable if there exist a measure
space $(\Omega,\Sigma, \mu)$ and operators $U_1 \in \cal L(E, L_p(\mu))$
and $U_2 \in \cal L(L_p(\mu), F^{**})$ such that $i_F \circ T= U_2 \circ U_1$
where $i_F: F \to F^{**}$ denotes the natural embedding.
\end{definition}

The $L_p$-factorable norm  of $T$ is defined by
 $\gamma_p(T):= \inf\{\Vert U_1 \Vert. \Vert U_2 \Vert\}$ where the infimum
 is  taken over all possible factorizations as above.

For detailed
discussion of $p$-summing operators, $(q,p)$-summing operators,
$p$-Schatten-von Neumann operators and $L_p$-factorable operators,
 we refer to \cite{DJT}, \cite{PIE} and \cite{TJ}.

 \medskip

We will now recall some basic facts on $C^*$-algebras and von-Neumann
algebras.
Let $\cal A$  be a $C^*$ algebra, we denote by $\cal A_h$ the set of
Hermitian (self adjoint) elements of $\cal A$. For $x \in \cal A$ and
$f \in \cal A^*$,  $xf$ (resp. $fx$) denotes the element of $\cal A^*$
defined  by $xf(y)=f(yx)$ (resp. $fx(y)=f(xy)$) for every $y \in \cal A$.

\begin{definition}
A von-Neumann algebra is said to be $\sigma$-finite if it admits at most
countably many orthogonal projections.
\end{definition}
We refer to \cite{KR} and \cite{TAK}
for some characterizations and examples of $\sigma$-finite von-Neumann
algebras.

\section{Main Theorem}

\begin{theorem}
Let $\cal A$ be a $C^*$-algebra and $T\in  \Pi_1 (\cal A, \ell_2)$.
 Then for every
$\varepsilon > 0$, there exists a Hilbert space $H$ and operators $J : A
\longrightarrow H$ and $K : H \longrightarrow \ell_2$
such that:
\begin{itemize}
\item[(1)]  $ T = K \circ J$;
\item[(2)] $\Vert J \Vert \leq 1$;
\item[(3)]  $K \in S_4(H, \ell_2)$ with $\sigma_4(K) \leq 2(1+\varepsilon)
\pi_1(T)$.
\end{itemize}
\end{theorem}

To prove this theorem, we will consider first the following particular case:
\begin{proposition}
Let $ \M$ be a $\sigma$-finite von-Neumann algebra,  $T: \M
\longrightarrow \ell_2$ be a weak* to weakly continuous absolutely summing
operator and $\varepsilon >0$.  Then there exist a Hilbert space $H$, operators
$J:  \M \longrightarrow H $ and $K: H \longrightarrow \ell_2$
 such  that:
 \begin{itemize}
\item[(1)] $T = K \circ J $;
\item[(2)] $\Vert J\Vert \leq 1$;
\item[(3)]
 $K \in  S_4 (H , \ell_2)$ with  $\sigma_4(K) \leq {2}
 (1+\varepsilon) \pi_1(T)$.
\end{itemize}
\end{proposition}

\begin{proof}
The proof is based on  the factorization technique used in \cite{Ran7}.  We
will repeat the argument for completeness.

Let $T \in  \Pi_1 ( \M, \ell_2)$ and assume that $T$ is weak*
to weakly continuous. Fix $\delta >0$ such that $(1+\delta)^{1/2}\leq
(1+\varepsilon)$.

 By \cite{Ran7} (Proposition 1.1) and \cite{Pis3}
 (Lemma 4.1),
there exists a normal positive functional $g$ on $\M$ such that
 $\Vert g \Vert \leq 1$ and
$$
\Vert T x \Vert \leq \pi_1 (T) g (\vert x \vert)\ \ \text{for every}\  x
\in  \M_h.
$$
Since the von-Neumann algebra $\M$ is $\sigma$-finite, there exists a
faithful normal functional  $f_0$ in $\cal M_*$
(see \cite{TAK} Proposition~II-3.19).
We can choose $f_0$ such that $\Vert f_0 \Vert \leq \delta$.
 Let $f = (g + f_0)/(1 +\delta)$; clearly
$\Vert f \Vert \leq 1$ and
$$
\Vert T x \Vert \leq (1+\delta)\pi_1 (T) f (\vert x \vert)\ \ \text{for every}\   x
\in  \cal M_h.
$$
From Lemma 2 of \cite{Ran7}, we deduce that
$$
\Vert T x  \Vert \leq 2(1+\delta) \pi_1 (T) \Vert x f\ +\ f x \Vert_{\cal M_*} \
  \ \text{for every}\  x
\in  \cal M.$$
As in \cite{Ran7}, we equip $\M$ with the scalar product
$$
\langle x, y \rangle = f \left(\frac {xy^* + y^* x}{2}\right).$$
Since $f$ is faithful, $\cal M$ with $\langle \cdot , \cdot \rangle$
is  pre-Hilbertian.  We denote the completion of this space by
 $L_2( \M, f)$  ( or  simply  $L_2(f)$).
From \cite{Ran7}, we have the following factorization:
$$
\xymatrix{
\M \ar[r]^-{T}  \ar[d]_-{J}  &\ell_2  &\\
L_2(f)    \ar[r]^-\theta            &L_2(f)^* \ar[r]^-
{J^*} & \M_* \ar[ul]_L}
$$
where $\theta (Jx) =  \langle \cdot \ , \ J (x^*)\rangle$
for every $x \in \M$;
 $L(\frac{xf + fx}{2}) = Tx$
\  for every $x \in  \M$ and $J$ is the inclusion map
(one can easily check as in \cite{Ran7} that
 $J^*\circ \theta\circ J(x)= (xf +fx)/2$).

Set $H :  =  L_2 (f)$ and $K:=  L \circ J^*\circ \ \theta$.
Clearly  (1) and (2) are satisfied.

\noindent
 To prove (3), let us consider the adjoint maps:
$$
\xymatrix{
& &\ell_2 \ar[r]  \ar[dl]_-{L^*}    &\cal M_{*}    \\
&\cal M \ar[r]^-J                 &L_2(f) \ar[r]^-\theta   &L_2(f)^{*}
\ar[u]_-{J^*}
}
$$
The proposition will be deduced from the following lemma:
\begin{lemma}
For every $p \geq 1$,
$K^*   \in  \ \Pi_{2p,p}(\ell_2, H^*)$ with
$\pi_{2p,p}(K^*)\leq \pi_p(T)^{\frac{1}{2}}\Vert L \Vert^{\frac{1}{2}}$.
\end{lemma}
To see the lemma,
let $(z_n)_n$ be a sequence in $\ell_2$ such that
$$
\sup\left\{\left(\sum_{n=1}^\infty \ \vert \langle z_n, z^* \rangle
\vert^p\right)^\frac 1p; \ \Vert z^* \Vert \leq 1\right\} \ = C < \infty.
$$
Then
$$
\sup\left\{  \left(\sum_{n=1}^\infty \vert \langle L^* (z_n), \xi
\rangle \vert^p \right)^ \frac 1p ; \xi \in \M^*, \ \Vert \xi \Vert \leq 1 \right\}
 \leq \Vert L \Vert C.
$$
Similarly,
$$
\sup\left\{ \left(\sum_{n=1}^\infty \vert \langle (L^*(z_n))^*,
 \xi \rangle \vert^p \right)^\frac 1p;
 \xi \in \M^*,\ \Vert \xi \Vert \leq 1 \right\} \leq \Vert L \Vert C
$$
where $\left(L^*(z_n)\right)^*$ is the adjoint of
the operator $L^*(z_n)$ in $\cal M$
for every $n \in  \N$.
Since $\left((L^*(z_n))^* \right)_n$ is a sequence in $\M$, one can
apply $T$. The fact that  $T$ is $p$-summing implies that
\begin{align*}
\left(\sum_{n=1}^\infty \Vert T (L^* (z_n)^*) \Vert^p \right)^\frac 1p
 &\leq \pi_p (T)
\sup\left\{ \left(\sum_{n=1}^\infty \vert \langle L^* (z_n)^*, \xi \rangle \vert^p \right)^
\frac 1p; \Vert \xi \Vert \leq 1 \right\} \\
&\leq \pi_p (T) \Vert L \Vert C.
\end{align*}
But since $(z_n)_n$ is bounded (in fact it is bounded by $C$) we get that
$$
\left(\sum_{n=1}^\infty \vert \langle T( (L^*(z_n))^*), z_n \rangle \vert^p \right)^\frac
1p \leq \pi_p (T) \Vert L \Vert C^2.
$$
Now for each $n \in  \N$,
\begin{align*}
\langle T \left( (L^*(z_n))^*\right), z_n \rangle
 &=\langle L \circ J^* \circ \theta \circ
J \left(L^* (z_n)^* \right), z_n \rangle \\
&= \langle \theta \circ J \left(L^*(z_n)^* \right), J \circ L^*(z_n)
\rangle \\
&= \langle J\left(L^*(z_n)\right), J \left(L^*(z_n)\right) \rangle \\
&= \Vert J \left(L^*(z_n)\right) \Vert_{L_{2} (f)}^2.
\end{align*}
So
$$
\left(\sum_{n=1}^\infty \Vert J \circ L^* (z_n)\Vert^{2p}\right)^ \frac 1p \leq \pi_p
(T) \Vert L \Vert C^2.
$$
Hence
$$
\left(\sum_{n=1}^\infty \Vert K^* (z_n)\Vert^{2p}\right)^\frac 1{2p} \leq \pi_p (T)^
\frac 12 \Vert L \Vert^\frac 12 C
$$
which shows that $K^* \in  \Pi_{2p,p}(\ell_2, H^*)$ with
$\pi_{2p,p}(K^*) \leq \pi_p (T)^\frac 12 \Vert L \Vert^\frac 12$.
The lemma is proved.

To complete the proof of the proposition,
 we apply the above lemma for $p = 2$;
we get that
$\pi_{4,2} (K^*) \leq \pi_2 (T)^\frac 12 \Vert L \Vert^\frac 12 \ $.
We note also from the proof of Theorem~1 of
\cite{Ran7} that the set $\{ xf + fx;\ x \in \M\}$ is
norm dense in $\M_*$ so from the estimate
 $\Vert L(\frac{xf + fx}{2}) \Vert = \Vert T x \Vert \leq {2}(1+\delta)
  \pi_1 (T) \Vert xf + fx\Vert_{\M_*} $ for every $x \in  \M $,
we get that
 $$\Vert L(xf +fx)\Vert \leq 4(1+ \delta)\pi_1(T) \Vert xf +fx \Vert_{\M_*}
 \ \ \text{for every}\  x \in \M.$$
  We conclude that
$\Vert L \Vert \leq 4 (1+\delta)\pi_1 (T) $ and therefore
 $\pi_{4,2}(K^*) \leq \pi_2(T)^\frac 12  \
 {2}(1+\delta)^{\frac 12} \pi_1(T)^\frac 12 \leq {2}
 (1+\varepsilon) \pi_1(T)$.

\noindent
From a result of Mitjagin (which appeared for the first time in a paper of
Kwapie\'n  \cite{KW2}; see also \cite{DJT} Theorem~10.3 or
 \cite{TJ} Proposition~11.8), the
space $\Pi_{4, 2}(\ell_2, H^*)$ is isometrically isomorphic to $S_4(\ell_2,
H^*)$ so $\sigma_4 (K^*) \leq {2}(1+\varepsilon) \pi_1 (T)$
 and from Proposition 4.5 of \cite{DJT}
(p. 80), $K \in S_4(H, \ell_2)$ with
 $\sigma_4(K) = \sigma_4(K^*) \leq {2}(1+\varepsilon)\pi_1 (T)$.
  The proof of the proposition  is complete.
\end{proof}

\noindent
{\bf Proof of Theorem~1.}
Assume first that $\cal A$ is separable and $T \in \Pi_1(\cal A, \ell_2)$. The
space ${\cal A}^{**}$ is a von-Neumann algebra and
 $T^{**} \in \Pi_1({\cal A}^{**}, \ell_2)$.
  Let $i_{\cal A}: \cal A \to {\cal A}^{**}$ be the natural embedding and
 choose $(a_n)_n$ a countable dense subset of $\cal A$. If $\M$ is the
 von-Neumann algebra generated by $\{i_{\cal A}(a_n); n\geq 1 \}$, then
 $\M$ is $\sigma$-finite. Also if we denote by $I$ the inclusion of $\M$
 into ${\cal A}^{**}$, then $I$ is weak* to weak* continuous.
 From  Proposition~1, the operator $T^{**}\circ I$ factors through a
 Hilbert space operator $K$ that belongs to the class $S_4$ and so does
 $T= T^{**}\circ I \circ i_{\cal A}$. One can easily verify that this
 factorization satisfies  the conclusion of the theorem.

For the general case, we will use ultraproduct technique.  Let $(\cal
A_s)_{s \in S}$ be the collection of all separable $C^*$-subalgebras
of $\cal A$.  As a particular case of Theorem 3.3 of \cite{SM}
 (which is the
$C^*$-version of Proposition 6.2 of \cite{He}), there exists a subset
$\Lambda$ of $S$ and an ultrafilter $\cal U$ on $\Lambda$ such that $\cal
A$ is (completely) isometric to a subspace of $(\cal A_s)_\cal U$.
Inspecting the proof of \cite{SM}, one notices that
in our case  $\Lambda = S$.

Let $T : \cal A \longrightarrow \ell_2$ be a $1$-summing operator and
 $i_s: \cal A_s \longrightarrow \cal A$ be the inclusion map.
It is clear that $T \circ i_s \in  \Pi_1 (\cal A_s, \ell_2)$ with
$\pi_1 (T\circ i_s) \leq \pi_1(T)$.  From the proposition above, there
exists a Hilbert space $H_s$ such that the following diagram commutes:
$$
\xymatrix{
\cal{A}_s\ar[dr]_-{J_s} \ar[rrr]^-{T_\circ i_s}&&& \ell_2\\
&H_s \ar[urr]_-{K_s}&&
}
$$
with $\Vert J_s \Vert \leq 1$ and $\sigma_4 (K_s) \leq 2(1 +
\varepsilon) \pi_1 (T)$.

From this, one can verify that the following diagram commutes:
$$
\xymatrix{
(\cal{A}_s)_{\cal U} \ar[dr]_-{(J_s)_{\cal U}} \ar[rr]^-{(T_\circ i_s)_{\cal
U}}&& (\ell_2)_{\cal U}\\
&(H_s)_{\cal U} \ar[ur]_-{(K_s)_{\cal U}}&
}
$$

It is clear that $\Vert (J_s)_\cal U \Vert \leq 1$ and since $S_4$ is a
maximal ideal operator, we get that
$(K_s)_{\cal U} \in  S_4\left( (H_s)_\cal U, (\ell_2)_\cal U \right)$
with $\sigma_4\left((K_s)_\cal U \right) \leq
 \lim\limits_{s,\cal U} \sigma_4(K_s)
\leq 2(1+ \varepsilon) \pi_1 (T)$ (see \cite{He} Theorem 8.1).

Let $Q: (\ell_2)_\cal U \longrightarrow \ell_2$ defined by
$Q \left((y_s)_s \right)=\text{weak}-\lim\limits_{s, \cal U} y_s$
 and $I: \cal A \longrightarrow (\cal A_s)_\cal U$ be the isometric embedding.
  We claim that
$Q\circ (T \circ i_s)_{\cal U}\circ I = T$.

 To see this, notice that  for every
$x \in \cal A$, $I(x)_s = 0$ if $x \notin \cal A _s$ and
$I(x)_s = x$ if $x \in \cal A_s$.  So $(T\circ i_s)_{\cal U}(Ix) =
(y_s)_{s \in S}$ where $y_s = 0$ if $x \notin \cal A_s$ and
$y_s = Tx$ if $x \in \cal A_s$ and by the definition of $Q$ the claim
follows.

 We get the conclusion of the theorem by setting
  $J=(J_s)_{\cal U}\circ I$, $K=Q\circ (K_s)_{\cal U}$ and
   $H=(H_s)_{\cal U}$.

 \qed

\smallskip

For the next simple extension of Theorem~1,
 we refer to \cite{UP} for definitions and
examples of $JB^*$-triples and $JBW^*$-triples.

\begin{corollary}
If $\cal A$ is a $JB^*$-triple then every absolutely summing operator from
$\cal A$ into $\ell_2$ factors through an operator that belongs to the
$4$-Schatten-von Neumann class.
\end{corollary}

\begin{proof}
Let $T: \cal A \rightarrow \ell _2$ be absolutely summing operator.  The
space $\cal A ^{**}$ is a $JBW^*$-triple.  But every $JBW^*$-triple is (as
Banach space) isometric to a complemented subspace of a von-Neumann algebra
(see \cite{CI2}).
From Theorem~1, $T^{**}$ (and consequently $T$) factors through an
operator that belongs to the class $S_4$.
\end{proof}

\begin{remark}
We remark that Lemma 1 is valid for any weak* to weakly continuous
absolutely summing operator from a $\sigma$-finite von-Neuman algebra
into a general Banach space;  in particular, the adjoint of any such
operator belongs to the class ideal $\Pi_{2p, p}$ for every $p \geq 1$.
\end{remark}
\bigskip
The following finite dimensional examples show that one can not improve
Theorem~1 to the case of p-Schatten-von Neumann class for $p < 4$.  The type
of operators considered below were suggested to the author by Pe\l czy\'nski.

\medskip

For $n \geq 1$,  $B(\ell^n_2)$ \ (resp. $HS(\ell^n_2)$)  denotes
the space of $n \times n$ matrices with the usual operator norm (resp. the
Hilbert-Schmidt norm).

\medskip
Let $I_n: B (\ell^n_2)  \longrightarrow HS(\ell^n_2)$ be the identity
operator and set \  $\alpha_n = \pi_1 (I_n)$.

\begin{theorem}
For every $n \geq 1$, let $T_n = {I_n}/{\alpha_n}$.
There exists an absolute constant $\beta > 0$ (independent of $n$)
 such that if $H$ is a
Hilbert space, $J \in  \cal L \left(B(\ell^n_2), H \right)$
and $K \in \cal L \left(H, HS (\ell^n_2)\right)$ satisfying:
\begin{itemize}
\item[(i)] $\Vert J \Vert \leq 1$;
\item[(ii)] $T_n = K \circ J$.
\end{itemize}
Then for every $p \geq 2$,
 $\displaystyle{\sigma_p (K) \geq \beta n^\frac{4-p}{2p}}$.
\end{theorem}

For the proof of this theorem, we will recall few well-known facts about the
operator $I_n$.

\begin{proposition}
1) There exists a universal constant $c>0$ such that
$\displaystyle
{\alpha_n = \pi_1 (I_n) \leq cn}$  for every $ n\geq 1$;

2) There exists a universal constant $c^\prime > 0$ such that
$\displaystyle
{\gamma_1 (I_n) \geq c^\prime n^\frac32}$ for every $n\geq 1$.
\end{proposition}

To prove the theorem, let $H$ be a Hilbert space and $J$ and $K$ be
operators as in the statement.  Since $HS(\ell^n_2)$ is a
finite dimensional Hilbert space,
$K: H \rightarrow HS(\ell^n_2)$ is a Hilbert-Schmidt operator.
Similarly, the adjoint
$K^*:  HS(\ell^n_2) \longrightarrow H$ is also a Hilbert-Schmidt
operator.  One can choose a probability space $(\Omega, \Sigma, \lambda)$
such that:
$$
\xymatrix{
HS(\ell_2^n) \ar [d]_-V \ar [r]^-{K^*}&H\\
L_\infty(\lambda) \ar [r]^-{i_2}  &L_2(\lambda) \ar [u]_-U
}
$$
with $\Vert V \Vert = 1$ and $\Vert  U \Vert =\pi_2(K^*) = \pi_2 (K)$.
Taking the adjoints,
$$
\xymatrix{
H \ar [d]_-{U^*} \ar [r]^-K&HS(\ell_2^n)\\
L_2(\lambda) \ar [r]^-{i_2^*}  &L_1(\lambda) \ar [u]_{V^*}
}
$$
Hence the operator $T_n$ factors through $L_1(\lambda)$ as follows:
$$
\xymatrix{
B(\ell_2^n) \ar[dr]_-{U_1} \ar[rr]^-{T_n}&&HS(\ell_2^n)\\
&L_1(\lambda) \ar[ur]_-{V^*}&
}
$$
where $ U_1 = i^*_2 \circ U^* \circ J$.  From the definition of
$\gamma_1(T_n)$, we get the following estimate:
\begin{align*}
\gamma_1 (T_n) & \leq \Vert \ U_1 \Vert.  \Vert V^* \Vert\\
& \leq \Vert i^*_2 \Vert. \Vert  U^* \Vert. \Vert J \Vert. \Vert
V^* \Vert\\
& \leq \Vert  U^* \Vert = \pi_2 (K).
\end{align*}
From the above proposition,
$\frac{c^\prime n^\frac12}c \
 \leq \frac {c^\prime n^\frac32} {\alpha_n} \leq \pi_2 (K)$.

 If  we set $\beta : = \frac{c^\prime}c$, we get
  $ \sigma_2(K)=\pi_2(K)\geq \beta n^\frac12$
 and the theorem is proved for the case $p = 2$.

For $p > 2$, note that $B(\ell^n_2)$ and $HS(\ell^n_2)$ are of dimension
$n^2$ so we can assume without loss of generality
that $\text{dim}(H) = n^2$.
 Let $(s_i (K))_{1\leq i \leq n^2}$ be the singular numbers of $K$.
  It is well known that for every $q >0$,
  $\sigma_q (K) = \left(\sum\limits_{i=1}^{n^2} s_i(K)^q\right)^\frac1q.$
Using Holder's inequality, we get for every $p > 2$,
\begin{align*}
\sigma_2 (K)  & =
  \left(\sum\limits_{i=1}^{n^2} s_i(K)^2 \right)^\frac12\\
& \leq  \left(\sum\limits_{i=1}^{n^2} s_i (K)^p \right)^\frac1p  .
\left(\sum\limits_{i=1}^{n^2} 1 \right)^{(1-\frac2p) \frac12}\\
& =  \sigma_p (K)   n^{1- \frac2p}.\\
\end{align*}
Hence
$\beta n^\frac12 \leq \sigma_2 (K) \leq \sigma_p (K) n^{1-\frac2p}$ which
implies that $\sigma_p (K) \geq \beta n^{-\frac12 + \frac2p} = \beta
n^\frac{4-p}{2p}$. The proof of the theorem is complete. \qed

\medskip
The operator $T_n$ satisfies $\pi_1 (T_n) = 1$  but any factorization
through any Hilbert space operator has large p-Schatten-von Neumann norm
for $p < 4$. This  shows that the class $S_4$ in the statement
of Theorem 1 cannot be improved.
\smallskip

The results above lead us to the question of characterizing operators from
a $C^*$-algebra into $\ell_2$ that can be factored through Hilbert-Schmidt
operators.

\begin{theorem}
Let $\cal A$ be a $C^*$-algebra.  An operator $T:  \cal A \longrightarrow
\ell_2$ factors through a Hilbert-Schmidt operator if and only if it is
$L_1$-factorable.
\end{theorem}

\begin{proof}
If $T$ factors through a Hilbert-Schmidt operator then it is
$L_1$-factorable since Hilbert-Schmidt operators are $L_1$-factorable.

Conversely, assume that $T$ is $L_1$-factorable i.e. there exits a measure
space $(\Omega, \Sigma, \lambda)$,
operators $ U_1: \cal A \longrightarrow L_1 (\Omega, \Sigma, \lambda)$
 and $ U_2: L_1(\Omega, \Sigma, \lambda) \longrightarrow \ell_2$
 such that $T = U_2 \circ U_1$.
From Grothendieck's theorem $ U_2$ is $1$-summing.  Since $L_1(\Omega,
\Sigma, \lambda)$ is of cotype $2$, $ U_1$ factors through a Hilbert
space  (see \cite{Pis3}) which shows that $T$ factors through a
Hilbert-Schmidt operator.
\end{proof}

\section{Measure of compactness.}

In this section, we will provide  an application of the main theorem to
measure compactness of any absolutely summing operator from $C^*$-algebras
into  Hilbert spaces.

Let $L$ be a normed linear space with norm $\| \ \cdot \ \|$ and $A$ be a
totally bounded set in $L$.

For any given $\varepsilon > 0$, we set $N_{\varepsilon} (A) : =$ the
infimum of integers $m$ such that there exists of subsets $E_1, E_2, \dots,
E_m$ of $L$ whose diameters do not exceed $2 \varepsilon$ and whose union
contains $A$ i.e,
 $$\bigcup^n_{k=1} E_k \supseteq A \text{ \ and \ } diam (E_k) \leq
2 \varepsilon.$$

\begin{definition}
$H_\varepsilon (A) : = \log_2 N_\varepsilon (A)$ is called the
$\varepsilon$-capacity of the set $A$.
\end{definition}
This definition was introduced by Kolmogorov and Tihomirov (among other
related notions) in \cite{KT}.

Our main result in this section answers positively a question raised by
Pe\l czy\'nski and can be viewed as a quantitative version of Theorem~1 of
\cite{Ran7}.

\begin{theorem}
There exists an absolute constant $C$ such that if $\cal A$ is a
 $C^*$-algebra
and $T \in \Pi_1 (\cal A, \ell_2)$ with $\pi_1 (T) \leq 1$, then for every
$\varepsilon > 0$,
$$ H_\varepsilon \left( T (B_{\cal A}) \right) \leq \frac
C{\varepsilon^4}.$$
\end{theorem}
We will show that Theorem~4 is a consequence of the following result.

\begin{theorem}
Let $H$ be a separable Hilbert space and $S \in S_p (H, \ell_2)$, then
for every
$\varepsilon > 0$, $$H_\varepsilon \left( S (B_H) \right) \leq
\frac{\sigma_p (S)^p \cdot \rho (p)}{\varepsilon^p}$$
 where $\rho (p) =
\left( \frac{8^p}p + \int^{8^{-p}}_0 \ln (\frac 1t) dt + 1 \right)^p $.
\end{theorem}
The proof is based on a notion of entropy of operators  introduced by Pietsch
(see \cite{PIE} p.~168).

\begin{definition}
Let $E$ and $F$ be Banach spaces and $S \in \cal  L(E, F)$.  The
$n$-th (outer) entropy number $e_n(S)$ of the operator $S$ is the
minimum of $\delta > 0$ such that there exists a finite sequence
 $y_1, y_2, \dots, y_q \in F$ with
$q \leq 2^{n-1}$ and  $S(B_E) \subseteq \bigcup^q_{i=1}\{ y_i +
\delta B_F \}$.
\end{definition}

Clearly $e_{n+1}(S) \leq e_n(S)$ for every operator $S$ and every $n \in
\N$.

For  diagonal Hilbert space operators, the following proposition was
proved by Pietsch.
\begin{proposition}
(\cite{PIE} p.~174)  Let $S \in \cal L(\ell_2)$ such that $S \left(
{(\xi_n)}_n  \right) = {(\alpha_n \xi_n)}_{n \geq 1}$ and
${(\alpha_n)}_n \in c_0$. Then
 $$\left( \sum^\infty_{n = 1} e_n (S)^p \right)^{\frac 1p} \leq K_p
\left( \sum_{n=1}^\infty | \alpha_n |^p \right)^{\frac 1p}.$$
\end{proposition}

\noindent
{\bf Proof of  Theorem~5.}

 Let  $1<p < \infty$ and $S \in S_p (H, \ell_2)$.  The
operator $S$ admits an orthonormal representation
$$(*) \quad S = \sum_{n=1}^\infty \alpha_n ( \ \cdot \ , \ h_n) f_n \ ,$$
where $(h_n)$ and $(f_n)$ are orthonormal sequences in $H$ and $\ell_2$
respectively and ${(\alpha_n)}_n \in \ell_p$.  We can choose this
representation so that $0 \leq \alpha_{n+1} \leq \alpha_n$ for all
admissible indices.  Let ${(e_n)}_n$ be the unit vector basis of $\ell_2$.
Let $Z = \overline{\text{span}}  \{h_n; n \in \N\}$ in $H$.  Since
$S(B_Z) = S(B_H)$, we can assume without loss of generality that $H = Z$.

Let $I : \ell_2 \longrightarrow H$ defined by $Ie_n = h_n $
for every $n \in \N$ and $J: \ell_2 \longrightarrow \ell_2$
 so that $J(f_n) = e_n$ for every $n \in \N$.  Let $\widetilde S
: = J \circ S \circ I$.
 Clearly $\widetilde S \in S_p (\ell_2,\ell_2)$,  $I$ and
$J$ are isometries.

For every $x \in \ell_2$, we have
\begin{align*}
\widetilde S x &= \sum^{\infty}_{n=1} \alpha_n (Ix, h_n) J (f_n)\\
&= \sum^{\infty}_{n=1} \alpha_n (x, I^*h_n) e_n \\
&= \sum^{\infty}_{n=1} \alpha_n (x, e_n) e_n  .
\end{align*}
So for every $x = {(x_n)}_n \in \ell_2$, $Sx = {(\alpha_n x_n)}_{n \geq
1}$.  Hence $\widetilde S$ satisfies the assumption of the above
proposition and therefore
$$\left( \sum^\infty_{n=1}
 \left( e_n (\widetilde S) \right)^p \right)^{\frac 1p} \leq K_p
\left( \sum^\infty_{n=1} | \alpha_n |^p \right)^{\frac 1p}
\leq K_p \sigma_p (S).$$
For $\varepsilon > 0$, define $k (\varepsilon) : = \max \{ k : e_k
(\widetilde S) \geq \varepsilon \}$.  We have
$$\left( K_p \sigma_p (S) \right)^p \geq \sum^{\infty}_{n=1} \left( e_n
(\widetilde S) \right)^p \geq \sum^{k(\varepsilon)}_{n=1} \left( e_n
(\widetilde S) \right)^p \geq \varepsilon^p k (\varepsilon)$$
so $$k (\varepsilon) \leq \left( \frac{K_p \sigma_p (S)}{\varepsilon}
\right)^p .$$
From the definition of $k(\varepsilon)$,  $e_{k(\varepsilon) + 1}
(\widetilde S) \leq \varepsilon$ and the definition of the $n$-th
entropy of $\widetilde S$
implies that there exists $\delta \leq \varepsilon$ and $\{ y_1, y_2,
\dots, y_q \} \subseteq \ell_2$, with $q \leq 2^{k(\varepsilon)}$ so that
$\widetilde S (B_{\ell_2}) \subset \{ y_1, y_2, \dots, y_q \} + \delta
B_{\ell_2}$ \ i.e.,
the set  $\widetilde S (B_{\ell_2})$ can be covered by
$2^{k(\varepsilon)}$ balls  of radius $\delta \leq \varepsilon$ so
$N_{\varepsilon} \left( \widetilde S (B_{\ell_2}) \right) \leq
2^{k(\varepsilon)}$ and
$$H_{\varepsilon} \left( \widetilde S (B_{\ell_2}) \right)
 \leq k(\varepsilon)
\leq \frac{{(\sigma_p (S) K_p)}^p}{\varepsilon^p} .$$
Now since $J$ is an isometry, $H_{\varepsilon} \left( \widetilde S
(B_{\ell_2}) \right) = H_{\varepsilon} \left( S \circ I (B_{\ell_2}) \right)$;
also by the definition of $I$, $I(B_{\ell_2}) = B_H$ so
$$H_{\varepsilon} \left( S (B_H) \right) = H_{\varepsilon} \left(
\widetilde S (B_{\ell_2}) \right) \leq \frac{\sigma_p (S)^p
{K_p}^p}{\varepsilon^p}$$
and setting $\rho (p) = {K_p}^p$, the theorem is proved.

  The estimate on
$K_p$ can be found in Pietsch's book \cite{PIE}( p.~174). \qed

\medskip
\noindent
{\bf  Proof of Theorem 4.}

\noindent
If $\cal A$ is  a $C^*$-algebra and $T \in
\Pi_1 (\cal A, \ell_2)$ with $\pi_1 (T) \leq
1$, then one can deduce from Theorem 1 and Theorem~5 that
$H_{\varepsilon} \left( T(B_{\cal A}) \right) \leq \frac{3^4
\rho(4)}{\varepsilon^4}$. In fact, one can choose $H$, $J$ and $K$ (as in
Theorem~1) so that $\sigma_4(K) \leq {3}$ so from Theorem~5,
$H_{\varepsilon}\left(K(B_{H})\right)\leq \frac{3^4
\rho(4)}{\varepsilon^4}$ and since $\Vert J \Vert \leq 1$,
$H_{\varepsilon}\left(T(B_{\cal A})\right)\leq \frac{3^4
\rho(4)}{\varepsilon^4}$.
Hence  if we set  $C = 3^4\rho(4)$, the proof of the theorem is complete.
 \qed

 \medskip
\noindent
{\bf Acknowledgments.} Some of the results in this paper were obtained
during the author's visit to the Institut of Mathematics of the Polish
Academy of Sciences, Warsaw, Poland. The author wishes to thank Prof.
 A. Pe\l czy\'nski for arranging the visit and for
several very insightful discussions concerning this work.
The author is also very grateful to Prof. S. Kwapie\'n for
some fruitful suggestions.


\end{document}